\documentclass[a4paper,11pt] {amsart}
\usepackage{amssymb,latexsym,amsmath}
\usepackage{amsfonts,amsbsy,bm}
\usepackage{color}
\usepackage{xcolor}
\usepackage{graphicx}
\usepackage{booktabs}
\usepackage{listings}
\newcommand{\bprop} {\begin{proposition}}
\newcommand{\eprop} {\end{proposition}}
\newcommand{\btheo} {\begin{theorem}}
\newcommand{\etheo} {\end{theorem}}
\newcommand{\blem} {\begin{lemma}}
\newcommand{\elem} {\end{lemma}}
\newcommand{\bcor} {\begin{corollary}}
\newcommand{\ecor} {\end{corollary}}

\newcommand{\Be}{\begin{equation}}
\newcommand{\Ee}{\end{equation}}
\newcommand{\Bea}{\begin{eqnarray}}
\newcommand{\Eea}{\end{eqnarray}}
\newcommand{\Bes}{\begin{equation*}}
\newcommand{\Ees}{\end{equation*}}
\newcommand{\Beas}{\begin{eqnarray*}}
\newcommand{\Eeas}{\end{eqnarray*}}
\newcommand{\Ba}{\begin{array}}
\newcommand{\Ea}{\end{array}}

\definecolor{Brown}{cmyk}{0,0.81,1,0.60}
\definecolor{OliveGreen}{cmyk}{0.64,0,0.95,0.40}
\definecolor{CadetBlue}{cmyk}{0.62,0.57,0.23,0}
\definecolor{lightlightgray}{gray}{0.91}

\begin{document}

\lstset{
	language=Python,                             
	basicstyle=\ttfamily,                   
	keywordstyle=\color{black},        
	commentstyle=\color{gray},              
	numbers=left,                           
	numberstyle=\tiny,                      
	stepnumber=1,                           
	numbersep=5pt,                          
	backgroundcolor=\color{lightlightgray}, 
	frame=none,                             
	tabsize=2,                              
	captionpos=b,                           
	breaklines=true,                        
	breakatwhitespace=false,                
	showspaces=false,                       
	showtabs=false,                         
	morekeywords={__global__, __device__},  
}

\theoremstyle{plain}
\newtheorem{theorem}{Theorem}[section]
\newtheorem{corollary}[theorem]{Corollary}
\newtheorem{lemma}[theorem]{Lemma}
\newtheorem{proposition}[theorem]{Proposition}

\theoremstyle{definition}
\newtheorem{definition}[theorem]{Definition}
\newtheorem{example}[theorem]{Example}

\theoremstyle{remark}
\newtheorem{remark}[theorem]{Remark}
\title[A two-parameter Gamma function
]{A two-parameter deformation of the gamma function, associated functions, and some related inequalities }
\author[A. Asare-Tuah]{Anton Asare-Tuah}
\address{Department of Mathematics, University of Ghana,\\ P. O. Box LG 62 Legon, Accra, Ghana}
\email{aasare-tuah@ug.edu.gh}
\author[E. Djabang]{Emmanuel Djabang}
\address{Department of Mathematics, University of Ghana,\\ P. O. Box LG 62 Legon, Accra, Ghana}
\email{EDjabang@ug.edu.gh}

\author[E. K. A. Schwinger]{Eyram K. A Schwinger}
\address{Department of Mathematics, University of Ghana,\\ P. O. Box LG 62 Legon, Accra, Ghana}
\email{eakschwinger@ug.edu.gh}
\author[B. F. Sehba]{Beno\^it Florent Sehba}
\address{Department of Mathematics, University of Ghana,\\ P. O. Box LG 62 Legon, Accra, Ghana}
\email{bfsehba@ug.edu.gh}
\author[R. A. Twum]{Ralph A. Twum}
\address{Department of Mathematics, University of Ghana,\\ P. O. Box LG 62 Legon, Accra, Ghana}
\email{ratwum@ug.edu.gh}


\maketitle

\begin{abstract}
In this paper, we introduce a new two-parameter deformation of the Gamma function that generalizes some existing Gamma-type functions in the literature. We study properties of this function that depend on the parameters.  We also prove some inequalities for the corresponding beta function and polygamma functions. 
\end{abstract}
\section{Introduction}
The aim of this paper is to introduce a new two-parameter gamma function, which is another deformation of the usual gamma function and related functions. We exhibit some partial differential equations (with respect to the parameter) satisfied by associated functions. We show some inequalities for the polygamma function and for ratios of the beta function. We compare our inequalities on the beta function to some existing inequalities in the literature.
\medskip

There are some deformations of the usual gamma function in the literature that have been intensively studied (see \cite{Diaz1,Diaz2,Djabang,Ege1,Ege2,Ege3,Gehlot2,Gehlot3,Gehlot4,Nantomah1,Wang}). Our proposed deformation generalizes two existing deformations in the literature. Let us recall these two deformations here:
\begin{itemize}
    \item The $k$-Gamma function of \cite{Diaz1} is defined by
    $$\Gamma_k(x)=\lim_{n\to\infty}\frac{n!k^n(nk)^{\frac{x}{k}-1}}{(x)_{n,k}},\quad k>0, x\in\mathbb{C}\setminus k\mathbb{Z}^-.$$
    \item The $\nu$-analogue of the Gamma function was introduced and studied in \cite{Djabang} and is defined by
    $$\Gamma^\nu(x)=\lim_{n\to\infty}\frac{n!\nu^n\left(\frac{n}{\nu}\right)^{\frac{x}{\nu}-1}}{(x)_{n,\nu}},\quad k>0, x\in\mathbb{C}\setminus \nu\mathbb{Z}^-.$$
\end{itemize}
Here, $(x)_{n,k}=x(x+k)\cdots (x+(n-1)k)$ is the so-called Pochhammer $k$-symbol.
Deformations of the classical Gamma function have applications in quantum mechanics and physics, combinatorics, special functions, and orthogonal polynomilas (see, for example, \cite{Deligne,Diaz1,Diaz3,Diaz4,Diaz2} and the references therein).
\medskip

In this note, we introduce the $(k,\nu)$-Gamma function as follows:
\Be\label{eq:gammanew}
\Gamma_{k,\nu}(x)=\lim_{n\to\infty}\frac{n!(k\nu)^n\left(\frac{nk}{\nu}\right)^{\frac{x}{k\nu}-1}}{(x)_{n,k\nu}},\quad k>0, \nu>0, x\in\mathbb{C}\setminus k\nu\mathbb{Z}^-.
\Ee
One clearly has that
$$\Gamma_{1,\nu}=\Gamma^\nu\quad\text{and}\quad \Gamma_{k,1}=\Gamma_k.$$
Note that the usual gamma function $\Gamma$ corresponds to $k=\nu=1$.
\medskip

In this work, we pay attention to some equations satisfied by the logarithm of the $(k,\nu)$ gamma function that depend on the parameters and then do not follow from the relation with the usual gamma function. We also prove some inequalities for the associated beta function and derivatives of the digamma function. 
\medskip

In the next section, we prove that the $(k,\nu)$ gamma function is unique in some sense by providing a Bohr-Mollerup-type theorem. In Section 3, we provide properties for the corresponding digamma function and some inequalities satisfied by the polygamma function. In Section 4, we obtain some inequalities for the beta function. We compare these inequalities to some others in the literature.
\medskip

The notation $\Gamma$, $B$, $\Psi$ ($k=\nu=1$), will always represent the usual gamma, beta and digamma functions, respectively. For more on these classical functions, we refer the reader to \cite{Abramowitz}.

\section{The $(k,\nu)$ Gamma function}
This section is more about basic facts on the deformed gamma function, we also prove here a Bohr-Mollerup’s theorem for the $(k,\nu)$ gamma function.
\medskip

We recall that for $x\in \mathbb{C}$, $a\in \mathbb{R}$ and $n\in\mathbb{N}$, the Pochhammer a-symbol is defined as 
\Be\label{eq:pochhammer}
(x)_{n,a}=x(x+a)\cdots(x+(n-1)a).
\Ee
We also recall our definition of $\Gamma_{k,\nu}$.
\begin{definition}
 For $k>0$ and $\nu>0$, the $(k,\nu)$ Gamma function $\Gamma_{k,\nu}$ is defined as   
 $$\Gamma_{k,\nu}(x)=\lim_{n\to\infty}\frac{n!(k\nu)^n\left(\frac{nk}{\nu}\right)^{\frac{x}{k\nu}-1}}{(x)_{n,k\nu}},\quad x\in\mathbb{C}\setminus k\nu\mathbb{Z}^-.$$
\end{definition}
We have the following basic facts.
\bprop
Let $k,l,\nu,\mu>0,$. For $x\in\mathbb{C}\setminus (k\nu\mathbb{Z}^-\cap l\mu\mathbb{Z}^-)$, and $n\in\mathbb{N}$, the following hold.
\begin{itemize}
    \item[(a)] $(x)_{n,l\mu}=\left(\frac{l\mu}{k\nu}\right)^n\left(\frac{k\nu x}{l\mu}\right)_{n,k\nu}$.
    \item[(b)] $\Gamma_{l,\mu}(x)=\left(\frac{l\nu}{k\mu}\right)^{\frac{x}{l\mu}-1}\Gamma_{k,\nu}\left(\frac{k\nu x}{l\mu}\right)$. 
    \item[(c)] $\Gamma_{k,\nu}(k\nu)=1$.
    \item[(d)] \Be\label{eq:gammmatogamma}\Gamma_{k,\nu}(x)=\left(\frac{k}{\nu}\right)^{\frac{x}{k\nu}-1}\Gamma\left(\frac{x}{k\nu}\right).\Ee
 
\end{itemize}
\eprop
We note that the function $(k,\nu)$ gamma is logarithmically convex. In effect, for $0<p,q<1$ with $p+q=1$, and $x,y>0$,
\Be\label{eq:holdergamma}
\Gamma_{k,\nu}(px+qy)\leq \left(\Gamma_{k,\nu}(x)\right)^p\left(\Gamma_{k,\nu}(y)\right)^q.
\Ee
Let us prove the following extension of Bohr-Mollerup's theorem (for the the functions $\Gamma(x)$ and $\Gamma_{k,1}(x)$, see \cite{Bohr,Diaz1}).
\btheo\label{thm:bohrgene}
Let $f:(0,\infty)\to (0,\infty)$, and let $k,\nu>0$. Suppose that
\begin{itemize}
    \item[(i)] $f(k\nu)=1$;
    \item[(ii)] $f(x+k\nu)=x\nu^{-2}f(x)$;
    \item[(iii)] $f$ is logarithmically convex. 
\end{itemize}
Then $f(x)=\Gamma_{k,\nu}(x)\quad\text{for all}\quad x\in (0,\infty)$.
\etheo
\begin{proof}
 We have by the definition of $\Gamma_{k,\nu}(x)$ that
 $$f(x)=\Gamma_{k,\nu}(x)\quad\text{if and only if}\quad \lim_{n\to\infty}\frac{(x)_{n,k\nu}f(x)}{n!(k\nu)^n\left(\frac{nk}{\nu}\right)^{\frac{x}{k\nu}-1}}=1.$$
 Fom the hypotheses on $f(x)$ and properties of $\Gamma_{k,\nu}(x)$, we can suppose that $0<x<k\nu$. 
 \medskip

 As $f$ is logarithmically convex, the following inequalities are satisfied:
 $$\frac{1}{k\nu}\ln\left(\frac{f(nk\nu+k\nu)}{f(nk\nu)}\right)\leq \frac{1}{x}\ln\left(\frac{f(nk\nu+k\nu+x)}{f(nk\nu+k\nu)}\right)\leq \frac{1}{k\nu}\ln\left(\frac{f(nk\nu+2k\nu)}{f(nk\nu+k\nu)}\right).$$
 As $f(x+k\nu)=x\nu^{-2}f(x)$, the above reduces to
 $$\frac{x}{k\nu}\ln\left(\frac{nk}{\nu}\right)\leq\ln\left(\frac{(x+nk\nu)\left(x+(n-1)k\nu\right)\cdots xf(x)}{n!\left(k\nu\right)^n}\right)\leq \frac{x}{k\nu}\ln\left(\frac{(n+1)k}{\nu}\right).$$
 That is
 $$0\leq \ln\left(\frac{(x+nk\nu)\left(x+(n-1)k\nu\right)\cdots xf(x)}{n!\left(k\nu\right)^n\left(\frac{nk}{\nu}\right)^{\frac{x}{k\nu}}}\right)\leq \ln\left(\left(\frac{n+1}{n}\right)^{\frac{x}{k\nu}}\right).$$
 As $\displaystyle\lim_{n\to\infty}\ln\left(\left(\frac{n+1}{n}\right)^{\frac{x}{k\nu}}\right)=0$, we conclude that
 $$\lim_{n\to\infty}\ln\left(\frac{(x+nk\nu)\left(x+(n-1)k\nu\right)\cdots xf(x)}{n!\left(k\nu\right)^n\left(\frac{nk}{\nu}\right)^{\frac{x}{k\nu}}}\right)=0.$$
 Hence
 $$\lim_{n\to\infty}\frac{(x)_{n,k\nu}f(x)}{n!(k\nu)^n\left(\frac{nk}{\nu}\right)^{\frac{x}{k\nu}-1}}=1.$$
 The proof is complete.
\end{proof}

\section{The $(k,\nu)$ Psi function and its derivatives}
In this section, we examine the $(k,\nu)$ psi function and its derivatives. 
\subsection{The $(k,\nu)$ Psi function and its properties}
\begin{definition}
Let $k,\nu>0$. The function $(k,\nu)$ psi (digamma) denoted $\Psi_{k,\nu}(x)$ is defined as the derivative of $\ln\left(\Gamma_{k,\nu}(x)\right)$.    
\end{definition}
We have the following that is needed in our next result.
\bprop\label{pro:psiexpression}
Let $k,\nu>0$. For any $x>0$, we have
\Be\label{eq:psiexpression}
\Psi_{k,\nu}(x)=\frac{1}{k\nu}\left(\ln k-\ln\nu-\gamma\right)-\frac{1}{x}-\sum_{j=1}^\infty\left(\frac{1}{x+jk\nu}-\frac{1}{jk\nu}\right).
\Ee

\eprop
\begin{proof}
    The identity follows from the definition of $\Psi_{k,\nu}(x)$ and the identity
    $$\frac{1}{\Gamma_{k,\nu}(x)}=\nu^{\frac{x}{k\nu}-1}k^{-\frac{x}{k\nu}}\frac{x}{\nu}e^{\frac{\gamma x}{k\nu}}\prod_{j=1}^\infty\left(1+\frac{x}{jk\nu}\right)e^{-\frac{x}{jk\nu}}.$$
\end{proof}
Let us prove the following about the anti-derivative of the $(k,\nu)$ digamma function.
\btheo\label{thm:pdepsi}
The function $\Psi(k,\nu,x)=\ln\left(\Gamma_{k,\nu}(x)\right)$ is a solution of the following nonlinear partial differential equations
\Be\label{eq:pdepsi1}
k^2\partial_k^2\Psi+2k\partial_k\Psi-x^2\partial_x^2\Psi=-1-\frac{x}{k\nu}
\Ee
and 
\Be\label{eq:pdepsi2}
\nu^2\partial_\nu^2\Psi+2\nu\partial_\nu\Psi-x^2\partial_x^2\Psi=2+\frac{x}{k\nu}.
\Ee
\etheo
\begin{proof}
Recall that 
$$\Psi(k,\nu,x)=\frac{x}{k\nu}\left(-\ln\nu+\ln k-\gamma\right)+2\ln\nu-\ln x-\sum_{j=1}^\infty\left(\ln\left(1+\frac{x}{jk\nu}\right)-\frac{x}{jk\nu}\right).$$
It follows that
$$\partial_k\Psi=\frac{x}{k^2}\left(\frac{\ln\nu-\ln k+\gamma+1}{\nu}+\sum_{j=1}^\infty\left(\frac{k}{x+jk\nu}-\frac{1}{j\nu}\right)\right);$$
$$\partial_k^2\Psi=-2\frac{x}{k^3}\left(\frac{\ln\nu-\ln k+\gamma+1}{\nu}+\sum_{j=1}^\infty\left(\frac{k}{x+jk\nu}-\frac{1}{j\nu}\right)\right)+\frac{x}{k^2}\left(-\frac{1}{k\nu}+\sum_{j=1}^\infty\frac{x}{(x+jk\nu)^2}\right).$$
Also,
$$\partial_x\Psi=\frac{1}{k\nu}\left(\ln k-\ln\nu-\gamma\right)-\frac{1}{x}-\sum_{j=1}^\infty\left(\frac{1}{x+jk\nu}-\frac{1}{jk\nu}\right)$$
and
\Be\label{eq:partial2psi}\partial_x^2\Psi=\sum_{j=0}^\infty\frac{1}{(x+jk\nu)^2}.\Ee
One easily get that
$$k^3\partial_k^2\Psi+2k^2\partial_k\Psi=kx\left(-\frac{1}{k\nu}+\sum_{j=1}^\infty\frac{x}{(x+jk\nu)^2}\right)$$
and so
$$k^3\partial_k^2\Psi+2k^2\partial_k\Psi-kx^2\partial_x^2\Psi=-k-\frac{x}{\nu}.$$
That is $\Psi$ satisfies (\ref{eq:pdepsi1}).

The proof for the second equation follows in the same way.
\end{proof}
\begin{remark}
\begin{itemize}
\item Equation (\ref{eq:partial2psi}) also tells us that $\Gamma_{k,\nu}(x)$ is logarithmically convex. We even derive from it that
$$\Gamma_{k,\nu}''(x)\Gamma_{k,\nu}(x)-\left(\Gamma_{k,\nu}'(x)\right)^2>0.$$
From the same (\ref{eq:partial2psi}), we deduce that $\Psi_{k,\nu}(x)$ is strictly increasing.
\item Obviously, the above result cannot be derived from a corresponding result for the case $k=\nu=1$.
\end{itemize}
\end{remark}

Let us prove the following identity.
\bprop
Let $k,\nu>0$. Assume that $x>0$.Then for any  $n\in\mathbb{N}\cup\{0\}$,
\Be\label{eq:diffpsi}
\Psi_{k,\nu}(x+(n+1)k\nu)-\Psi_{k,\nu}(x)=\sum_{j=0}^n\frac{1}{x+jk\nu}.
\Ee
Moreover,
\Be\label{eq:difflim}\lim_{n\to\infty}\left(\Psi_{k,\nu}(x+(n+1)k\nu)-\frac{1}{k\nu}\ln n\right)=\frac{1}{k\nu}\ln\left(\frac{k}{\nu}\right).
\Ee
\eprop
\begin{proof}
From (\ref{eq:psiexpression}) and (\ref{eq:diffpsi}), we obtain 
\Beas
\Psi_{k,\nu}(x+(n+1)k\nu) &=& \Psi_{k,\nu}(x)+\sum_{j=0}^n\frac{1}{x+jk\nu}\\ &=& \frac{1}{k\nu}\left(\ln\left(\frac{k}{\nu}\right)-\gamma\right)-\sum_{j=n+1}^\infty\left(\frac{1}{x+jk\nu}-\frac{1}{jk\nu}\right)+\sum_{j=1}^n\frac{1}{jk\nu}.
\Eeas
Hence
$$\Psi_{k,\nu}(x+(n+1)k\nu)-\frac{1}{k\nu}\ln n=\frac{1}{k\nu}\left(\ln\left(\frac{k}{\nu}\right)-\gamma\right)-\sum_{j=n+1}^\infty\left(\frac{1}{x+jk\nu}-\frac{1}{jk\nu}\right)+\frac{1}{k\nu}\left(\sum_{j=1}^n\frac{1}{j}-\ln n\right).$$
Taking the limit as $n\to\infty$, it easily follows that
$$\lim_{n\to\infty}\left(\Psi_{k,\nu}(x+(n+1)k\nu)-\frac{1}{k\nu}\ln n\right)=\frac{1}{k\nu}\ln\left(\frac{k}{\nu}\right).$$
The proof is complete.
\end{proof}
\subsection{$(k,\nu)$ Polygamma function and some inequalities}
We start this part with a definition.
\begin{definition}
Let $m>0$ be an integer. The  $(k,\nu)$ m-polygamma function $\Psi_{k,\nu}^{(m)}$ is the $m^{th}$ derivative of the $(k,\nu)$ Psi function i.e.
$$\Psi_{k,\nu}^{(m)}(x)=\frac{d^m}{dx^m}\Psi_{k,\nu}^{(m)}(x).$$
\end{definition}
From Proposition \ref{pro:psiexpression}, we obtain the following.
\bprop\label{prop:polygaidi}
Let $k,\nu>0$, and let $m>0$ be an integer. Then for $x>0$, we have
\Beas
\Psi_{k,\nu}^{(m)}(x) &=& (-1)^{m+1}m!\sum_{n=0}^\infty\frac{1}{(x+nk\nu)^{m+1}}\\ &=& (-1)^{m+1}\int_0^\infty\frac{t^me^{-xt}}{1-e^{-k\nu t}}dt.
\Eeas
\eprop
From the above, we make the following useful observation.
\begin{equation} \label{eq:polygavar}
\Psi_{k,\nu}^{(m)}\quad\text{is}\, 
\begin{cases}
    \text{increasing and concave} & \text{if m is even}\\
    \text{decreasing and convex} & \text{if m is odd}
\end{cases}
\end{equation}
We have the following.
\bprop\label{prop:polygamineq1}
Let $k,\nu>0$, and let $m\in\mathbb{N}\cup\{0\}$. Put $\Psi_{k,\nu}^{(0)}=\Psi_{k,\nu}$ and $\Psi_{k,\nu}^{(-1)}=\ln\Gamma_{k,\nu}$. Then for $0<x<y$, we have the following
\Be\label{eq:polygamineq11}
\frac{\Psi_{k,\nu}^{(2m-1)}(y)-\Psi_{k,\nu}^{(2m-1)}(x)}{y-x}
\leq \Psi_{k,\nu}^{(2m)}\left(\frac{x+y}{2}\right)\Ee
and 
\Be\label{eq:polygamineq12}
\Psi_{k,\nu}^{(2m+1)}\left(\frac{x+y}{2}\right)\leq \frac{\Psi_{k,\nu}^{(2m)}(y)-\Psi_{k,\nu}^{(2m)}(x)}{y-x}
.
\Ee
\eprop
\begin{proof}
To prove (\ref{eq:polygamineq11}), we use the fact that $\Psi_{k,\nu}^{(2m)}$ is concave and the Jensen's inequality as follows
\Beas
\Psi_{k,\nu}^{(2m)}\left(\frac{x+y}{2}\right) &\geq& \frac{1}{y-x}\int_x^y\Psi_{k,\nu}^{(2m)}(t)dt\\ &=& \frac{\Psi_{k,\nu}^{(2m-1)}(y)-\Psi_{k,\nu}^{(2m-1)}(x)}{y-x}.
\Eeas
The inequality (\ref{eq:polygamineq12}) follows with the same type of arguments.
\end{proof}
Let us prove the following.
\bprop\label{prop:polygamineq2}
Let $k,\nu>0$, and let $m\in\mathbb{N}\cup\{0\}$. Put $\Psi_{k,\nu}^{(0)}=\Psi_{k,\nu}$ and $\Psi_{k,\nu}^{(-1)}=\ln\Gamma_{k,\nu}$. Then for $0<x<y$, we have the following 
\Be\label{eq:polygamineq21}
\frac{\Psi_{k,\nu}^{(2m)}(x)+\Psi_{k,\nu}^{(2m)}(y)}{2}\leq \frac{\Psi_{k,\nu}^{(2m-1)}(y)-\Psi_{k,\nu}^{(2m-1)}(x)}{y-x}
\Ee
and
\Be\label{eq:polygamineq22}
\frac{\Psi_{k,\nu}^{(2m)}(y)-\Psi_{k,\nu}^{(2m)}(x)}{y-x}\leq\frac{\Psi_{k,\nu}^{(2m+1)}(x)+\Psi_{k,\nu}^{(2m+1)}(y)}{2}.
\Ee
\eprop
\begin{proof}
We only prove the first inequality because the second one uses the same argument. For this, we consider the function
$$F(y)=\Psi_{k,\nu}^{(2m-1)}(y)-\Psi_{k,\nu}^{(2m-1)}(x)-(y-x)\frac{\Psi_{k,\nu}^{(2m)}(x)+\Psi_{k,\nu}^{(2m)}(y)}{2}.$$
Clearly,
\Beas
F'(y) &=& \Psi_{k,\nu}^{(2m)}(y)-\frac{\Psi_{k,\nu}^{(2m)}(x)+\Psi_{k,\nu}^{(2m)}(y)}{2}-(y-x)\frac{\Psi_{k,\nu}^{(2m+1)}(y)}{2}\\ &=& \frac{y-x}{2}\left[\frac{\Psi_{k,\nu}^{(2m)}(y)-\Psi_{k,\nu}^{(2m)}(x)}{y-x}-\Psi_{k,\nu}^{(2m+1)}(y)\right]\\ &=& \frac{y-x}{2}\left(\Psi_{k,\nu}^{(2m+1)}(c)-\Psi_{k,\nu}^{(2m+1)}(y)\right)
\Eeas
for some $c\in(x,y)$, thanks to the Mean Value Theorem. But by (\ref{eq:polygavar}), $\Psi_{k,\nu}^{(2m+1)}$ is decreasing. Hence $F'(y)\geq 0$, that is $F$ is an increasing function. As $x<y$, we deduce that $F(y)\geq F(x)=0$. That is
$$\Psi_{k,\nu}^{(2m)}(y)-\Psi_{k,\nu}^{(2m)}(x)-(y-x)\frac{\Psi_{k,\nu}^{(2m-1)}(x)+\Psi_{k,\nu}^{(2m-1)}(y)}{2}\geq 0$$
or equivalently,
$$\frac{\Psi_{k,\nu}^{(2m)}(x)+\Psi_{k,\nu}^{(2m)}(y)}{2}\leq \frac{\Psi_{k,\nu}^{(2m-1)}(y)-\Psi_{k,\nu}^{(2m-1)}(x)}{y-x}.$$
\end{proof}
Let us also prove the following.
\bprop\label{prop:polygamineq3}
Let $k,\nu>0$, $\theta\geq 0$, and let $m\in\mathbb{N}\cup\{0\}$. Then for $x,y>0$, and $r>1$, we have the following
\Be\label{eq:polygamineq31}
\frac{\left[\Psi_{k,\nu}^{(2m)}(\theta+x)\right]^r}{\Psi_{k,\nu}^{(2m)}(\theta+rx)}\leq \frac{\left[\Psi_{k,\nu}^{(2m)}(\theta+x+y)\right]^r}{\Psi_{k,\nu}^{(2m)}(\theta+r(x+y))}
\Ee
and 
\Be\label{eq:polygamineq32}
\frac{\left[\Psi_{k,\nu}^{(2m+1)}(\theta+x+y)\right]^r}{\Psi_{k,\nu}^{(2m+1)}(\theta+r(x+y))}\leq \frac{\left[\Psi_{k,\nu}^{(2m+1)}(\theta+x)\right]^r}{\Psi_{k,\nu}^{(2m+1)}(\theta+rx)}.
\Ee
\eprop
\begin{proof}
We only prove the first inequality as the second one follows in the same way.
\medskip

Let $$F(t)=\ln\left(\frac{\left(\Psi_{k,\nu}^{(2m)}(\theta+t)\right)^r}{\Psi_{k,\nu}^{(2m)}(\theta+rt)}\right)=\ln\left(\Psi_{k,\nu}^{(2m)}(\theta+t)\right)^r-\ln\left(\Psi_{k,\nu}^{(2m)}(\theta+rt)\right),\quad t>0.$$
Then
\Beas F'(t) &=& r\frac{\Psi_{k,\nu}^{(2m+1)}(\theta+t)}{\Psi_{k,\nu}^{(2m)}(\theta+t)}-r\frac{\Psi_{k,\nu}^{(2m+1)}(\theta+rt)}{\Psi_{k,\nu}^{(2m)}(\theta+rt)}\\ &=& r\frac{\Psi_{k,\nu}^{(2m+1)}(\theta+t)\Psi_{k,\nu}^{(2m)}(\theta+rt)-\Psi_{k,\nu}^{(2m)}(\theta+t)\Psi_{k,\nu}^{(2m+1)}(\theta+rt)}{\Psi_{k,\nu}^{(2m)}(\theta+t)\Psi_{k,\nu}^{(2m)}(\theta+rt)}.
\Eeas
By (\ref{eq:polygavar}), we have that
$$\Psi_{k,\nu}^{(2m)}(\theta+t)\leq \Psi_{k,\nu}^{(2m)}(\theta+rt)$$
and $$\Psi_{k,\nu}^{(2m+1)}(\theta+rt)\leq \Psi_{k,\nu}^{(2m+1)}(\theta+t).$$
Hence $F'(t)\geq 0$. That is, $t\mapsto F(t)$ is increasing and so $t\mapsto \frac{\left(\Psi_{k,\nu}^{(2m)}(\theta+t)\right)^r}{\Psi_{k,\nu}^{(2m)}(\theta+rt)}$ is increasing. Thus,
$$\frac{\left[\Psi_{k,\nu}^{(2m)}(\theta+x)\right]^r}{\Psi_{k,\nu}^{(2m)}(\theta+rx)}\leq \frac{\left[\Psi_{k,\nu}^{(2m)}(\theta+x+y)\right]^r}{\Psi_{k,\nu}^{(2m)}(\theta+r(x+y))}.$$
\end{proof}
The following is obtained as above.
\bprop\label{prop:polygamineq4}
Let $k,\nu>0$, $\theta\geq 0$, and let $m\in\mathbb{N}\cup\{0\}$. Then for $x,y>0$, and $r<1$, we have the following
\Be\label{eq:polygamineq41}
\frac{\left[\Psi_{k,\nu}^{(2m)}(\theta+x+y)\right]^r}{\Psi_{k,\nu}^{(2m)}(\theta+r(x+y))}\leq \frac{\left[\Psi_{k,\nu}^{(2m)}(\theta+x)\right]^r}{\Psi_{k,\nu}^{(2m)}(\theta+rx)}
\Ee
and 
\Be\label{eq:polygamineq42}
\frac{\left[\Psi_{k,\nu}^{(2m+1)}(\theta+x)\right]^r}{\Psi_{k,\nu}^{(2m+1)}(\theta+rx)}\leq \frac{\left[\Psi_{k,\nu}^{(2m+1)}(\theta+x+y)\right]^r}{\Psi_{k,\nu}^{(2m+1)}(\theta+r(x+y))}.
\Ee
\eprop
\section{Some inequalities for the beta function}
We prove some inequalities for the ratios of the beta function.
\begin{definition}
 Let $k,\nu>0$. The $(k,\nu)$ Beta function $B_{k,\nu}$ is defined by
 \Be\label{eq:knubetadef}
 B_{k,\nu}(x,y)=\frac{\Gamma_{k,\nu}(x)\Gamma_{k,\nu}(y)}{\Gamma_{k,\nu}(x+y)}.
 \Ee
\end{definition}
Using (\ref{eq:gammmatogamma}), one obtains that $$B_{k,\nu}(x,y)=\frac{\nu}{k}B(x,y).$$
\subsection{Inequalities for the ratios of the beta function}
Let us prove the following two-sided inequality.
\btheo\label{thm:twoside1}
Let $k,\nu>0$. For $x_1,x_2,y>0$ with $x_1<x_2$, the following holds.
\Be\label{eq:twoside1}
\frac{x_2+y}{x_1+y}\left(\frac{x_1}{x_2}\right)^{\frac{y}{k\nu}+1}<\frac{B_{k,\nu}(x_2,y)}{B_{k,\nu}(x_1,y)}<\frac{x_2+y}{x_1+y}\left(\frac{x_1}{x_2}\right)\left(\frac{x_1+y+k\nu}{x_2+y+k\nu}\right)^{\frac{y}{k\nu}}.
\Ee
\etheo
\begin{proof}
By putting $a=\frac{x}{k\nu}$  and $b=\frac{y}{k\nu}$,
we can reduce the matter to the usual Beta function $B(a,b)$. We then only have to prove that for $0<a_1<a_2$ and $b>0$, we have
\Be\label{eq:twoside1beta}
\frac{a_2+b}{a_1+b}\left(\frac{a_1}{a_2}\right)^{b+1}<\frac{B(a_2,b)}{B(a_1,b)}<\frac{a_2+b}{a_1+b}\left(\frac{a_1}{a_2}\right)\left(\frac{a_1+b+1}{a_2+b+1}\right)^{b}.
\Ee
Recall that $$B(a,b)=\int_0^1t^{a-1}(1-t)^{b-1}dt.$$
Put $g(a,b)=B(a+1,b)$. Then
$$\frac{\partial g}{\partial a}=\int_0^1t^{a}(1-t)^{b-1}\ln(t)dt.$$
Using that for $0<t<1$,
$$\frac{t-1}{t}<\ln t<t-1,$$
we obtain
$$-\int_0^1t^{a-1}(1-t)^{b}dt< \frac{\partial g}{\partial a}<-\int_0^1t^{a}(1-t)^{b}dt.$$
That is
$$-B(a,b+1)< \frac{\partial g}{\partial a}<-B(a+1,b+1)$$
or equivalently,
$$-\frac{b}{a}g(a,b)< \frac{\partial g}{\partial a}<-\frac{b}{a+b+1}g(a,b).$$
That is
$$-\frac{b}{a}<\frac{\frac{\partial g}{\partial a}}{g(a,b)}<-\frac{b}{a+b+1}.$$
Integrating from $a=a_1$ to $a=a_2$, we obtain after simplification, that
$$\ln\left(\frac{a_1}{a_2}\right)^b<\ln\left(\frac{B(a_2+1,b)}{B(a_1+1,b)}\right)<\ln\left(\frac{a_1+b+1}{a_2+b+1}\right)^b,$$
that is 
$$\left(\frac{a_1}{a_2}\right)^b<\frac{B(a_2+1,b)}{B(a_1+1,b)}<\left(\frac{a_1+b+1}{a_2+b+1}\right)^b.$$
Using that $B(x+1,y)=\frac{x}{x+y}B(x,y)$, we obtain that the above is equivalent to (\ref{eq:twoside1beta}).
The proof is complete.
\end{proof}
\begin{remark}
Note that following \cite[Theorem 3.13]{From}, we also have the following upper estimate
\Be\label{eq:twoside2}
\frac{B_{k,\nu}(x_2,y)}{B_{k,\nu}(x_1,y)}<\left(\frac{x_1+y}{x_2+y}\right)^{\frac{y}{k\nu}},\quad y>0,\,\, 0<x_1<x_2.
\Ee 
\end{remark}
The following generalizes \cite[Theorem 3.12 and Theorem 3.16]{From}. As the proof follows the lines of \cite{From}, it is left to the reader.
\btheo
Let $k,\nu>0$. For $y>0$ and $0<x_1<x_2$, the following hold.
\Be\label{eq:twoside31}
\left(\frac{x_2}{k\nu}\right)^{\frac{x_2}{k\nu}-1}\left(\frac{x_1}{k\nu}\right)^{1-\frac{x_1}{k\nu}}\left(\frac{x_2+y}{k\nu}\right)^{1-\frac{x_2+y}{k\nu}}\left(\frac{x_1+y}{k\nu}\right)^{\frac{x_1+y}{k\nu}-1}\leq \frac{B_{k,\nu}(x_2,y)}{B_{k,\nu}(x_1,y)}
\Ee
and

\Be\label{eq:twoside32}
\frac{B_{k,\nu}(x_2,y)}{B_{k,\nu}(x_1,y)}\leq \left(\frac{x_2}{k\nu}\right)^{\frac{x_2}{k\nu}}\left(\frac{x_1}{k\nu}\right)^{-\frac{x_1}{k\nu}}\left(\frac{x_2+y}{k\nu}\right)^{-\frac{x_2+y}{k\nu}}\left(\frac{x_1+y}{k\nu}\right)^{\frac{x_1+y}{k\nu}}.
\Ee
\etheo
\begin{remark}
\begin{itemize}
\item[(a)]Following \cite[Remark 4]{From}, the upper bound in (\ref{eq:twoside32}) is better than the upper bound in (\ref{eq:twoside2}).  
\item[(b)] We have that the upper bound in (\ref{eq:twoside1}) is better that the one in (\ref{eq:twoside2}).
\begin{proof}
    Let us show that
    $$\frac{b+y}{a+y}\left(\frac{a}{b}\right)\left(\frac{a+y+1}{b+y+1}\right)^{y}< \left(\frac{a+y}{b+y}\right)^{y}$$
    or equivalently,
    $$\left(\frac{a}{b}\right)\left(\frac{a+y+1}{b+y+1}\right)^{y}< \left(\frac{a+y}{b+y}\right)^{y+1}.$$
    Put $f(y)=\left(\frac{a}{b}\right)\left(\frac{a+y+1}{b+y+1}\right)^{y}$ and $g(y)=\left(\frac{a+y}{b+y}\right)^{y+1}$. Then we have 
    \begin{eqnarray*}
        \ln\left(f(y)\right)-\ln\left(g(y)\right) &=& \ln a-\ln b+y\left(\ln(a+y+1)-\ln(b+y+1)\right)-(y+1)\left(\ln(b+y)-\ln(a+y)\right)\\ &=& -\int_a^b\frac{1}{t}dt-y\int_{a+y+1}^{b+y+1}\frac{1}{t}dt+(y+1)\int_{a+y}^{b+y}\frac{1}{t}dt\\ &=& \int_a^b\left(-\frac{1}{t}+\frac{y+1}{t+y}-\frac{y}{t+y+1}\right)dt\\ &=& -\int_a^b\frac{y^2+y}{t(t+y)(t+y+1)}dt\\ &<& 0.
    \end{eqnarray*}
    Hence $f(y)<g(y)$ which is enough to conclude that the upper bound in (\ref{eq:twoside1}) is better than the one in (\ref{eq:twoside2}).   
\end{proof}
\item[(c)] The lower bound in (\ref{eq:twoside31}) is better than the one in (\ref{eq:twoside1}).
\begin{proof}
It is enough to show that for $y>0$, $0<a<b<\infty$, $$b^{b-1}a^{1-a}(b+y)^{1-b-y}(a+y)^{a+y-1}>\frac{b+y}{a+y}\left(\frac{a}{b}\right)^{y+1}.$$ 
Put $$f(y)=b^{b-1}a^{1-a}(b+y)^{1-b-y}(a+y)^{a+y-1}$$ and $$g(y)=\frac{b+y}{a+y}\left(\frac{a}{b}\right)^{y+1}.$$
One easily checks that
\begin{eqnarray*}
h(y) &:=& \ln f(y)-\ln g(y)=(a+y)\ln\left(1+\frac{y}{a}\right)-(b+y)\ln\left(1+\frac{y}{b}\right).    
\end{eqnarray*}
As the function $t\mapsto (t+y)\ln\left(1+\frac{y}{t}\right)$ is strictly deacrasing, we conclude that $h(y)>0$ for any $y>0$. Thus, $f(y)>g(y)$.
\end{proof}

\item[(d)] However, the inequality between the upper bound in (\ref{eq:twoside32}) and the one in (\ref{eq:twoside1}) depends on the entries.

To compare the two bounds, we only have to compare the following functions of $(a,b,y)$:
$$A(a,b,y)=\frac{b^{b+1}}{a^{a+1}}\frac{(a+y)^{a+y+1}}{(b+y)^{b+y+1}}$$
and $$B(a,b,y)=\left(\frac{a+y+1}{b+y+1}\right)^{y}.$$
We note that if $A(a,b,y)-B(a,b,y)>0$, then the upper bound in (\ref{eq:twoside32}) is better than the one in (\ref{eq:twoside1}). If $A(a,b,y)-B(a,b,y)<0$, then the upper bound in (\ref{eq:twoside1}) is better than the one in (\ref{eq:twoside32}).
\medskip

Put $f(x)=(x+1)\ln(x+1)-(x+y+1)\ln(x+y)+y\ln(x+y+1)$. Observe that $\ln(A(a,b,y))-\ln(B(a,b,y))=f(b)-f(a)$.

We have $$f'(x)=\ln(x+1)-\ln(x+y)+\frac{y(x+y)-(x+y+1)}{(x+y)(x+y+1)}.$$
We note that $$-1\leq \frac{y(x+y)-(x+y+1)}{(x+y)(x+y+1)}\leq 1\quad\forall x,y>0.$$
Hence, we have for example that
\begin{itemize}
    \item for $0<x\leq 1$ and $y\geq 2e-1$, $f$ is decreasing. This implies that for $0<a<b\leq 1$ and $y\geq 2e-1$, $A(a,b,y)<B(a,b,y)$.
    \item For $0<y<e^{-1}$, and $0<x<\frac{1-ey}{e-1}$, $f$ is increasing. Hence, for $0<y<e^{-1}$, and $0<a<b<\frac{1-ey}{e-1}$, $A(a,b,y)>B(a,b,y)$.
\end{itemize}
\end{itemize}
\end{remark}
\medskip
\section{Compliance with Ethical Standards}
\begin{itemize}
\item {\bf Funding}

There is no funding support to declare.

\item {\bf Disclosure of potential conflicts of interest}

The authors have no relevant financial or non-financial interests to disclose.

\item {\bf Author Contributions}

All authors contributed to the conception and design of the study.  All authors read and approved the final manuscript.

\item {\bf Data availability statements}

Data sharing is not applicable to this article, as no data was created or analyzed in this study. 
\end{itemize}

\bibliographystyle{plain}

\end{document}